\documentclass[11 pt]{article}
\usepackage{amsmath}

\usepackage{newtxtext,newtxmath}

\usepackage{newtxmath,newtxtext}
\usepackage{microtype}

\usepackage[dvipsnames]{xcolor}

\usepackage[all]{xy}
\usepackage{graphicx,xspace,bm,bbm,tikz,pgfplots}

\usetikzlibrary{intersections,calc,patterns,angles,quotes}
\usetikzlibrary{decorations.shapes,arrows,calc,decorations.markings,shapes,decorations.pathreplacing,shapes.arrows,arrows.meta	}
\usepackage{tikz-3dplot}

\usepackage{url}

\usepackage{subfig}
\usepackage[skip=1pt,font=small]{caption}

\usepackage{xkeyval}
\usepackage{fullpage}

\newcommand{\Span}[1]{\operatorname{span} #1}

\newcommand{\diag}{\operatorname{diag}}
\newcommand{\R}{\mathbbm{R}}

\newcommand{\vx}{{\bf x}}
\newcommand{\va}{{\bf a}}
\newcommand{\vb}{{\bf b}}
\newcommand{\vv}{{\bf v}}
\newcommand{\vu}{{\bf u}}

\newcounter{para}

\newcommand {\bCE}{\left\lbrace \begin{aligned}}
\newcommand {\eCE}{ \end{aligned}\right.}

\usetikzlibrary{external}

\definecolor{darkred}{rgb}{0.9,0.1,0.1}

\newcommand{\nvar}[2]{%
    \newlength{#1}
    \setlength{#1}{#2}
}

\nvar{\dg}{0.3cm}

\nvar{\ddxc}{1.5cm}

\newtheorem{theorem}{Theorem}[section]

\date{}
\title{A Homotopy Method for Motion Planning.} 
\author{Shenyu Liu and Mohamed Ali Belabbas$^*$\thanks{$^{*}$Shenyu Liu  and Mohamed Ali Belabbas are with the department of Electrical and Computer Engineering and the Coordinated Science Laboratory, University of Illinois, Urbana-Champaign.
        {\tt\small sliu113,belabbas@illinois.edu}}}
\begin{document}
 \maketitle
\begin{abstract}
We propose a novel method for motion planning and illustrate its implementation on several canonical examples. The core novel idea underlying the method is to define a metric for which a path of minimal length is an admissible path, that is path that respects the various constraints imposed by the environment and the physics of the system on its dynamics. To be more precise, our method takes as input a control system with holonomic and non-holonomic constraints, an initial and final point in configuration space, a description of obstacles to avoid, and an initial trajectory for the system, called a sketch.  This initial trajectory does not need to meet the constraints, except for the obstacle avoidance constraints. The constraints are then encoded in an inner product, which is used to  deform (via a homotopy) the initial sketch  into an admissible trajectory from which controls realizing the transfer can be obtained. We illustrate the method on various examples, including vehicle motion with obstacles and a two-link manipulator problem.
\end{abstract}

\section{Introduction}

A fundamental problem in robotic motion planning is to find a trajectory which meets the various constraints stemming from the system's dynamics, which can be of holonomic or non-holonomic type,  and obstacle avoidance constraints, which include constraints on the magnitude of some of the variables describing the system (e.g., a maximal turning radius), or  obstacles present in physical space.
We propose here a new method to find a trajectory which takes into account all the above constraints--we call such a trajectory {\it admissible}--and illustrate its performance on several examples. The method is a homotopy method: given an initial state and a final desired state, $\vx_i$ and $\vx_f$ respectively, and an arbitrary curve joining $\vx_i$ to $\vx_f$ in state-space, the method {\it deforms} the curve into an {admissible} curve joining $\vx_i$ to $\vx_f$. We presented a preliminary version of this method, with only  non-holonomic constraints, in~\cite{Belabbas2017NewMF}. In this paper, we restrict the presentation to systems {\it affine in the control}, and leave the general case to subsequent work.
We also refer the readers to the website\footnote{\url{https://publish.illinois.edu/belabbas/motion-planning/}} for slides, sample Matlab code and examples showcasing the method.

The problem of motion planning in robotics and control is a canonical problem, and many methods have been proposed over the years. For this reason, we can only give here a very partial overview of the current state of the field, and emphasize that the method we propose is built on a rather different set of ideas. A large subset of the methods is focused on non-holonomic dynamics, since this problem is by itself difficult and with a long history~\cite{laumond1998robot,latombe2012robot, choset2005principles, Lav06}. Many of the proposed methods are based on the  use of sinusoidal driving signals; the basic  relation underlying these methods is the system approximation \begin{equation*}\dot x= \lim_{\omega \to \infty} \left( \sqrt{\omega}   \sin (\omega t) f_1(x)+ \sqrt{\omega} \cos(\omega t) f_2(x)\right) \Leftrightarrow \dot x  = [f_1,f_2](x),\end{equation*} where $[f_1,f_2]$ is the Lie bracket~\cite{do1992riemannian} of the vector fields $f_1,f_2$. Indeed, this insight is at the basis of the work of Brockett~\cite{brockett1989rectification}, Murray et al.~\cite{murray1994mathematical}, Laferriere and Sussman~\cite{lafferriere1993differential}.  Furthermore, interesting recent work shows that some special functions--which can be thought as generalizations of harmonic functions---play a distinguished role in solving under-actuated control problems~\cite{ gauthier2014minimal}. 

For control and verification of hybrid systems in general, we refer to~\cite{tomlin2003computational} and for a recent survey of motion planning for self-driving vehicles in urban environment, we refer to~\cite{paden2016survey}.  %
 Other approaches of interest to obtain  feasible trajectories for given problems and dynamics including random sampling-based~\cite{doi:10.1177/0278364911406761}   graph-based~\cite{1241712}, and optimization-based approaches~\cite{7041375} and approaches based on solvers for nonlinear dynamics.

\section{Background and problem set-up}

We present some  background and notation needed to explain the method. We refer to as vehicle/robot/plant  whose motion we desire to plan as {\it the system}. The system is assumed to obey the controlled dynamics 
\begin{equation}\label{eq:maindyn}\dot {\bf x} = \sum_{i=1}^p u_i f_i(\vx),\end{equation} 
where $\vx \in M$ with $M$ a (at least locally) differentiable manifold called the {\it configuration space}, $f_i(\vx)$ the {\it actuation vector fields} and $\vu:=(u_1,\ldots,u_p) \in \R^p$ the controls. We refer to as {\it workspace} the physical environment in which the system lives. We denote by $\Span_x \{ g_i \}$ the (real) vector space spanned by the vectors $g_i(x)$.

We call a {\it curve} in configuration space a piecewise differentiable function $\vx(t):[0,T] \to M$, where $T>0$, and refer to $\vx(0)$ and $\vx(T)$  as start-point and end-point, respectively, of $\vx(t)$. We refer to  them collectively as {\it end-points}. We call the {\it image} of a curve a {\it path}; a path is thus a geometric object (a collection of ''contiguous states'') and the times at which each point in a path is visited are not specified.

 A {\it fixed end-points homotopy} between the two curves $\vx_1(t)$ and $\vx_2(t)$ with the same end-points (i.e., $\vx_1(0)=\vx_2(0)$ and $\vx_1(T)=\vx_2(T)$) is a differentiable function $\vv(s,t):[0,\infty)\times[0,T] \to M$ with the properties:
\begin{align*}
\vv(s,0) &= \vx_1(0)&\mbox{for all } s \geq 0\\
\vv(s,T) &= \vx_1(T)&\mbox{for all } s \geq 0	
\end{align*}

The {\it length} of a curve $\vx(t)$ is defined with respect to an norm on the tangent bundle $TM$ of $M$. In the following, one can assume that $M=\R^n$ and the tangent space of $M$ at $\vx \in M$, denoted by $T_{\vx}M$ is also $\R^n$. A {\it Riemannian inner product} on $M$ is an given by piecewise differentiable symmetric positive definite bilinear form $G(\vx):T_{\vx}M \times T_{\vx}M \to \R.$. With a slight abuse of notation, we also denote by $G(\vx)$ its matrix representation in coordinates. Hence, we can think of $G(\vx)$ as an $\vx$-dependent positive definite symmetric matrix.

The {\it length} of a curve $p(t)$ is then given by \begin{equation}\label{eq:dfL}L(\vx):= \int_0^T \sqrt{\dot \vx^\top(t) G(\vx(t)) \dot \vx(t)} dt.\end{equation}

Finally, we introduce the Christoffels' symbols associated to $G(\vx)$. To this end, denote by $g_{ij}$ the $ij$th entry of the matrix representation $G(\vx)$, and by $g^{ij}$ the $ij$th entry of the matrix $G^{-1}(\vx)$. The Christoffel's symbols are

\begin{equation}\label{eq:defchristoffel} 
\Gamma_{jk}^i(\vx) := \frac{1}{2}\sum_l  g^{il}\left( \frac{\partial g_{lj}}{\partial x_k} +  \frac{\partial g_{lk}}{\partial x_j} -   \frac{\partial g_{jk}}{\partial x_l}\right)
 \end{equation}

\paragraph{Problem definition.} The problem that the method MotionSketch solves is the following: given a configuration space $M$, a set of holonomic, non-holonomic and obstacle avoidance constraints, an initial state $\vx_i$ and a desired final state $\vx_f$, provide a curve $\vx(t):[0,T] \to M$ which respects these constraints and so that $\vx(0)=\vx_i$, and $\vx(T)=\vx_f$, and provide the control $\vu$ that drive a control system from $\vx_i$ to $\vx_f$. From now on, we  normalize the time $T$ to be equal to one; this is done for simplicity of exposition, and all the results below are easily extended to the case of arbitrary $T$.
We recall that a curve that meets the constraints is an {\bf admissible curve}.

\begin{figure}
\begin{center}
\includegraphics[scale=1.2]{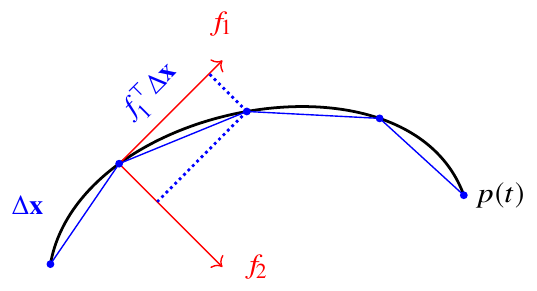}
\end{center}
\caption{Length of a discretized curve.}\label{fig:Lill}
\end{figure}

\paragraph{Length of a curve.} In order to provide an intuitive justification of the method, we first revisit the definition of the generalized length of a curve given a Riemannian metric in~\ref{eq:dfL}. See also Fig.~\ref{fig:Lill}. Since $G(\vx)$  is positive definite for all $\vx \in M$, we can factor it as $G(\vx)=F(\vx) D(\vx) F^\top(\vx)$, where $D(\vx)$ is a positive definite diagonal matrix, and $F(\vx)^\top F(\vx) = I$ (i.e., $F(\vx)$ is an orthogonal matrix.) Let $\vx(t):[0,1]\to M$ be a differentiable curve and let $0=t_0<t_1<\ldots<t_{l+1}=1$ provide subdivisions of the unit interval. We can then approximate $$\dot \vx(t_i) \simeq \frac{1}{\Delta t_i} (\vx(t_{i+1})-\vx(t_i))= \frac{1}{\Delta t_i} (\Delta \vx(t_i)),$$ where $\Delta t_i = t_{i+1}-t_i$, and the second equality defines $\Delta \vx(t_i)$. Using these relations, we can approximate the length of $\vx(t)$ as \begin{align*}L(\vx)&\simeq \sum_{i=1}^l \sqrt{(\frac{\Delta \vx(t_i)}{\Delta t_i} )^\top F(t_i)D(t_i)F(t_i)\frac{\Delta \vx(t_i)}{\Delta t_i}   } \Delta t_i \\
&\simeq \sum_{i=1}^l \sqrt{ (F(t_i)^\top\Delta \vx(t_i))^\top D(t_i)(F(t_i)\Delta \vx(t_i))}, 
\end{align*}
where we set $D(t_i):=D(\vx(t_i))$ and $F(t_i):=F(\vx(t_i))$. Since $F$ is an orthogonal matrix, we can think of  $F^\top \Delta \vx$ the vector of coordinates describing $\Delta \vx$ in the basis spanned by the column vectors of $F$; more precisely, if we set $f_k$ to be the $k$th column of $F$ and set $\Delta \vx_k(t_i) = f_k^\top \Delta \vx(t_i)$, then we have $\Delta \vx(t_i) = \sum_k f_k \Delta \vx_k(t_i)$. Now denote by $d_k^2$ the $k$th diagonal entry of $D$ (recall that $D$ has positive diagonal entries). We obtain 
$$L(\vx) \simeq \sum_i \sum_{k=1}^n \Delta \vx(t_i)_k d_k(t_i).$$ Hence, by adjusting the $d_i$ and the $f_k$ appropriately, we can {\it adjust which infinitesimal directions for a curve yield a larger length.} We show how this can be brought to bear on motion planning problems below. 

\section{The method MotionSketch}

The method contains the three following steps:
\begin{enumerate}
\item Encode the constraints of the motion planning problem (obstacles, holonomic, nonholonomic and dynamical constraints) into a Riemannian inner product.
\item Provide a curve in configuration space between the initial and final desired states. This curve, which we call the {\it sketch}, does {\it not} need to meet the holonomic, non-holonomic and dynamical constraints, but is required to avoid obstacles. Numerically solve the geometric heat flow (GHF), defined below, equation with the sketch as initial condition.
\item Extract the controls from the solution of the GHF.
\end{enumerate}

We now elaborate on the three items.

\subsection{Step 1: Encoding the constraints in a Riemannian inner product}

We start with holonomic/non-holonomic constraints.

\subsubsection{Holonomic and non-holonomic constraints}
Holonomic constraints can be formulated as a set of equations
\[
q_i(\vx)=0,\quad i=1,2,\cdots, m_h
\]
For each $i$ and an infinitesimally small motion $\delta \vx$, we have the approximation
$
q_i(\vx_0+\delta \vx)\approx q_i(\vx_0)+\frac{\partial q_i}{\partial \vx}\delta \vx.
$
In order to respect the constraint, $\delta \vx$ needs to satisfy $q_i(\vx_0+\delta \vx)=q_i(\vx_0)=0$, thus we have $\frac{\partial q_i}{\partial \vx}\delta \vx=0$. This means that for $\vx(t)$ to be an admissible curve, the direction of motion $\delta \vx$ needs to be orthogonal to the vectors $\frac{\partial q_i}{\partial \vx}$ for all $i$; in other words, it means the {\it undesirable}  directions of motion are $\Span \left\{\frac{\partial q_i}{\partial \vx}\right\}$.

We now turn our attention to non-holonomic constraints, which we assume are formulated  as a set of constraints on the allowed velocities $\dot \vx$ when at state $x$ as follows:
\[
\dot \vx^\top f_{c,j}(\vx)=0,\quad j=1,2,\cdots, m_n.
\] 
The non-holonomic character of the constraints, which is reflected in the fact that they cannot be expressed as $\frac{d}{dt} q_{n}(\vx) = 0$ for some function $q_{n}(\vx)$, does {\it not} play any particular function insofar our local encoding of the constraints is concerned; in fact,  the undesirable directions of motion are easily seen to be in this case  $\Span \left\{f_{c,j}(\vx)\right\}$.

Non-holonomic constraints can be presented as above, e.g. as non-slippage constraints, but they can also be encoded in the dynamics of the system, which is then called non-holononic. For this latter case, consider given the system of Eq.~\eqref{eq:maindyn}. We set $f_{f,i}=f_i$ and  $f_{c,j}$ to be  the $m_n$ vectors  orthogonal (for the Euclidean inner product) to $f_{f,i}$ for all $i=1,\cdots,p$.

\paragraph{Encoding the constraints} We set ${\bar p}:=n-m_n-m_h$. We define the $n\times(n-{\bar p})$ matrix $\bar F_c$ as the matrix with first $m_h$ columns given by $\frac{\partial q_i}{\partial \vx}$  and the next $m_n$ columns given by the $f_{c,j}$. We assume that $\bar F_c(\vx)$ is of constant rank almost everywhere in $M$, and we denote this rank by $l$, and set $p:=n-l$. If $m_h+m_n = l$, it is of {\it full column rank}, and we set $F_c(\vx):= \bar F_c(\vx)$. Otherwise $m_h+m_n> l$ and the constraints are not independent, in the sense that satisfying a {\it subset} of the constraints insures that {\it all} constraints are met. We set  $F_c(\vx)$ to be a $n \times {l}$ matrix whose column span equals the column span of $\bar F_c(\vx)$. Such matrix can be obtained, e.g., via the Gram-Schmidt process. Notice that $F_c$ is of full column rank $l=n-p$ and the {\it column space of $F_c$ contains all the undesirable directions of motion}.

Next, find a rank $p$ matrix $F_f(\vx)\in\R^{n\times p}$ such that
\[
F_f(\vx)^\top F_c(\vx)=0,
\]
which again can be found using the Gram-Schmidt process. The column space of $F_f(\vx)$ contains all the directions in which the system can move when at state $\vx$. Note that in the absence of holonomic constraints, we can start with defining $F_f$ with columns $f_i$ as in Eq.~\eqref{eq:maindyn} and choose $F_c$ the satisfy the above relation. Set
\begin{equation}\label{F(x)}
 F(\vx)=\begin{pmatrix}
|&|\\
F_c(\vx)&F_f(\vx)\\
|&|
\end{pmatrix}   
\end{equation}

Then $F(\vx)\in\R^{n\times n}$ and we define
\begin{equation}\label{H(x)}
H(\vx)=F(\vx)DF^\top (\vx)
\end{equation}
where $D=\diag([\underbrace{k \cdots k}_{n-p}\underbrace{ 1 \cdots 1}_p])$ is a constant matrix. Note that this $k$ is exactly the $d^2$ discussed in the Section II.b. In practice, we take $k$ to be of the order of $10\sim 1000$. 

Using the interpretation of the length functional given in the previous section, it is easy to see that if $\dot \vx$ is a  direction that respects the constraints, it is not multiplied by $k$ in  the inner product $\dot \vx^\top H(\vx)\dot \vx$ with $H$ defined via \eqref{H(x)}, so $\dot \vx^\top H(\vx)\dot \vx$ will not be scaled by $k$. On the other hand, if $\dot \vx$ is a  direction that violates a constraint,  it has some components lying in $\Span F_c(\vx)$,  and consequently $\dot \vx^\top H(\vx)\dot x$ is large.

Finally, we record here that the partial derivative of $H$ is given by
\[
\frac{\partial H}{\partial x_i}(\vx)=2F D\frac{\partial F^\top}{\partial x_i}(\vx),
\]
which is needed for the computation of the Christoffels symbols.
\subsubsection{Obstacle constraints}

We described obstacles  $\Omega_i\subset \R^n$ in configuration space via functions $r_i:M \to \R$ according to \[
\Omega_i:=\{\vx\in \R^n:r_i(\vx)\leq 0\}
\]
The boundary of an obstacle is thus  $\partial \Omega_i=\{\vx\in \R^n:r_i(\vx)=0\}$. We incorporate obstacles in the Riemannian inner product via a barrier function $b(\vx)=\sum_i b_i(\vx)$ with  the following properties:
\begin{enumerate}
\item Each $b_i(\vx)$ is positive and differentiable for all $\vx\in \R^n\backslash \Omega_i$ 
\item $b_i(\vx)\to\infty$ as $\vx\to\partial \Omega_i$,
\item $b(\vx)=1$ when $\vx$ is far away from all $\Omega_i$.
\end{enumerate}

\begin{figure}\centering
\subfloat[]{
\includegraphics[width=.51\columnwidth]{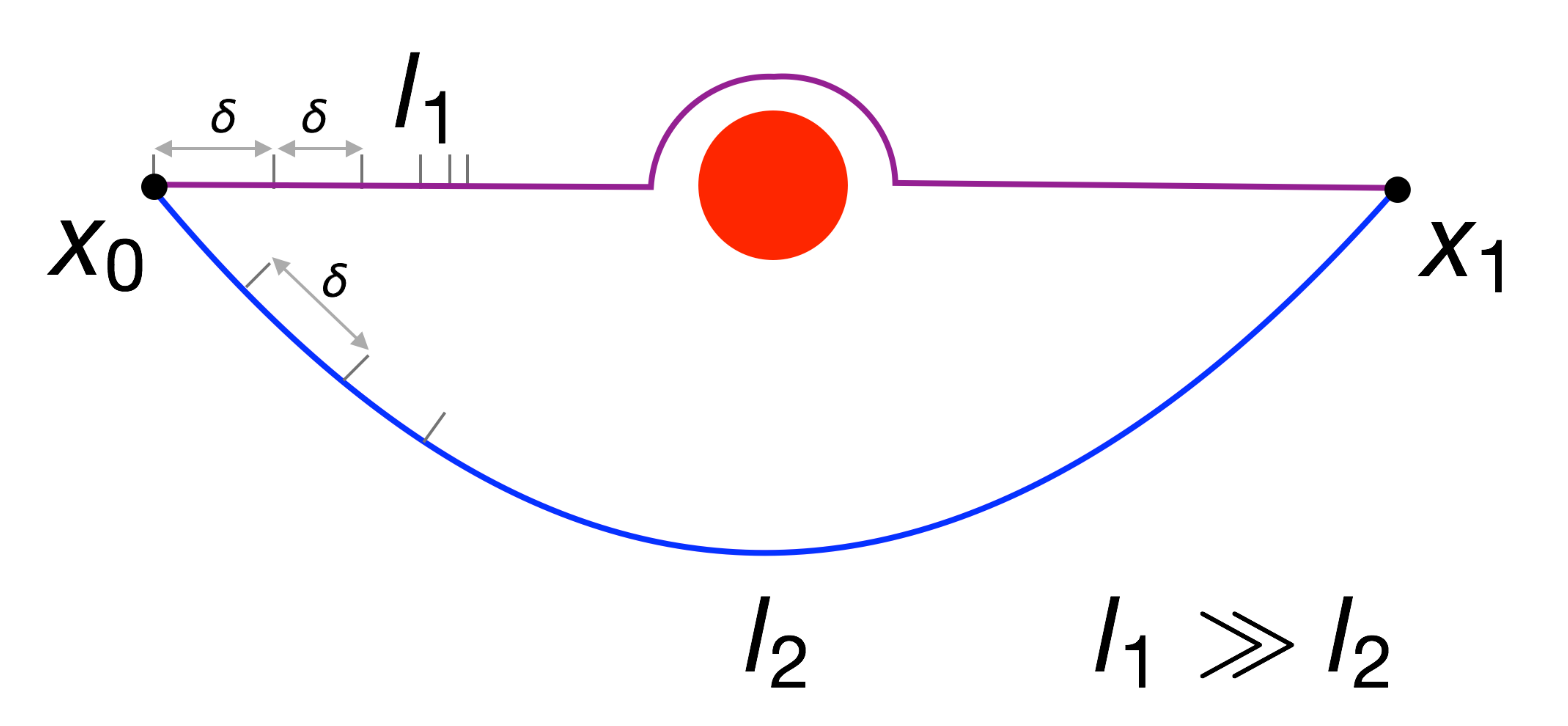}\label{fig:obstaclew}}
\subfloat[]{\includegraphics[scale=1.2]{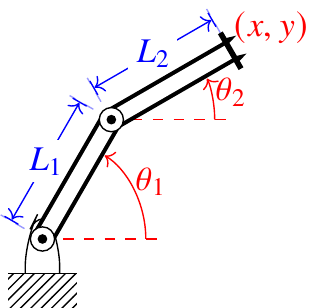}
\label{fig:2links}
}	
\caption{\small (a).The length $l_1$ of the path passing near the obstacle is much larger that the length $l_2$ of the path staying far from the obstacles when the metric is scaled with $b(\vx)$. (b) Two-links articulated arm can be described as a system with $4$ degrees of freedom and $2$ holonomic constraints relating the position $(x,y)$ of the tip to the joint angles $\theta_1,\theta_2$.}
\end{figure}

The idea is that we would like $b_i(\vx)$ to be large when $\vx$ is in the vicinity of $\Omega_i$, and becomes infinite if $\vx\in\partial\Omega_i$. Thus if we multiply the metric tensor by $b(\vx)$, the length of a path that is in the vicinity of an obstacle is much larger than the length of a path that steer well-clear of the obstacle, where quantifying ``well-clear'' is of course dependent on the choice of $b_i(\vx)$ and how quickly it decays near the boundary of the obstacle. We illustrate this in Fig.~\ref{fig:obstaclew}.

Such  functions $b_i$ are also known as \emph{barrier functions} in the optimization literature \cite{Nocedal99}. In the case when obstacles are balls, that is, $\Omega= \cup_{i=1}^l \{\vx\in \R^n:|\vx-c_i|\leq r_i\}$, one candidate of such $b(\vx)$ function will be a modification of penalty function from avoidance control \cite{Leitmann1980}:
\begin{equation}\label{b_ball}
b(\vx)=1+\sum_{i=1}^l\left(\min\left\{0,\frac{\vert \vx-c_i\vert^2-R_i^2}{\vert \vx-c_i\vert^2-r_i^2}\right\}\right)^2
\end{equation}
where $R_i$ is  such that  $r_i<R_i$ for all $i=1,2,\cdots,l$, and $R_i$ can be thought of as a {\it radius of detection} of the obstacle, in the sense that outside this radius, the obstacle does not affect the metric. Notice that $b(x)$ defined in \eqref{b_ball} satisfies the 3 properties mentioned earlier. The derivative of $b$ is also not hard to compute. Note that one can cover any obstacles with balls and use the above barrier function as a default approach.

\subsubsection{simultaneous multi-vehicle path planning}
Suppose there are $l$ vehicles and each of them has its own state $\vx_j=\begin{pmatrix}x_{1j},x_{2j},\cdots x_{nj}\end{pmatrix}^\top\in \mathbb R^n$ and the dynamics is $\dot\vx_j=F_j(\vx_j)\vu_j$. The $j$-th vehicle is supposed to drive from $\vx_j(0)=\va_j$ to $\vx_j(T)=\vb_j$. Denote $\vx^\top=\begin{pmatrix} \vx_1^\top&\cdots\vx_l^\top\end{pmatrix}$ and $\vu^\top=\begin{pmatrix} \vu_1^\top&\cdots\vu_l^\top\end{pmatrix}$, then the system of multi-vehicle has total dimension of $lm$ and initial and final states 
\[
\vx_i=\begin{pmatrix}\va_1\\\vdots\\\va_l\end{pmatrix},\quad\vx_f=\begin{pmatrix}\vb_1\\\vdots\\\vb_l\end{pmatrix}.
\]
and the overall dynamics is
\begin{equation}\label{bigF}
\dot \vx=\diag(F_1(\vx_1),\cdots,F_l(\vx_l))\vu:=F(\vx)\vu.
\end{equation}
While planning the path for all $l$ vehicles, they are also supposed to avoid collision with each other. In case of planar vehicles where $(x_{1,j},x_{2,j})$ represents the $xy$-coordinate of the $j$-th vehicle, collision between the $j,k$-th vehicles is avoided if 
\begin{equation}\label{collision_dist}
(x_{1j}-x_{1k})^2+(x_{2j}-x_{2k})^2\geq r_c^2,
\end{equation}
where $r_c$ is a safety radius guaranteeing collision-free between two vehicles. Thus the \eqref{b_ball}-like barrier function induced from \eqref{collision_dist} is
\[
b_c(\vx)=\sum_{j\neq k}\left(\min\left\{0,\frac{(x_{1j}-x_{1k})^2+(x_{2j}-x_{2k})^2-R^2}{(x_{1j}-x_{1k})^2+(x_{2j}-x_{2k})^2- r_c^2}\right\}\right)^2
\]
Thus, whenever two vehicles are too close ($(x_{1j}-x_{1k})^2+(x_{2j}-x_{2k})^2\leq R^2$), $b_c(\vx)$ becomes large and the metric at this state of vehicles is large. Notice that if we perform path planning for each individual vehicle first while treating the other vehicles as obstacles, the avoidance problem becomes dynamic in the sense that now the obstacles are moving with respect to time. Yet in our method avoidance of collision between vehicles and avoidance of static obstacles are processed in similar way and the result is promising as one can see later in our example.

In addition, Because $F(\vx)$ in \eqref{bigF} is block diagonal, $H$ defined via \eqref{H(x)} is also block diagonal and its $j$-th block only involves $\vx_j$. As a result, inverse of $H$ is in complexity of $O(lm^3)$ and computing $\frac{\partial H}{\partial x_i}$ for multi-vehicle has the same complexity as that for single vehicle. As a result, in each iteration of solving the numerical GHF equation, the complexity of computing all the Christoffel symbols is linear in $l$, the number of total vehicles.

\subsubsection{The inner product with three type of constraints}
We now formally define the inner product used in the method: given $H(\vx)$ as defined above from holonomic and non-holonomic constraints, and $b(\vx)$ a barrier function for the obstacles, we set
 
\[
G(\vx):=b(\vx)H(\vx)
\]
 With this construction, the partial derivatives of $G(\vx)$ can be computed using the chain rule:
$
\frac{\partial}{\partial x_i}G(\vx)=\frac{\partial b}{\partial x_i}(\vx)H(\vx)+b(\vx)\frac{\partial H}{\partial x_i}(\vx).
$
Hence the Christoffel symbols in \eqref{eq:defchristoffel} can be computed solely based on the values $H,\frac{\partial H}{\partial x_i},b,\frac{\partial b}{\partial x_i}$ at each state $\vx$.

\subsubsection{Examples}
\paragraph{The two-links manipulator}

In this example we consider a two-links manipulator in the plane, see Fig.~\ref{fig:2links}. The working space, in terms of the position of the tool tip $(x,y)$, is a subset of $\R^2$. The configuration space when the joint angles are also taken into account can be treated as a subset of $\R^4$. This system has 2 degrees of freedom and we can easily obtain the holonomic constraints:
\begin{equation}\label{eq:twoarmconst}
\left\{\begin{array}{c}
q_1(\vx)=L_1\cos(\theta_1)+L_2\cos(\theta_2)-x=0\\
q_2(\vx)=L_1\sin(\theta_1)+L_2\sin(\theta_2)-y=0
\end{array}
\right.
\end{equation}
Taking differential of the two  constraints, we find\begin{small}
\begin{equation*}
\frac{\partial q_1}{\partial \vx}=   (-1,\ 0,\ -L_1 \sin \theta_1,\-L_2 \sin \theta_2)^\top,\\ \frac{\partial q_2}{\partial \vx}=(  0,\ -1,\ L_1 \cos \theta_1 ,\ L_2 \cos \theta_2 )^\top
\end{equation*}
\end{small}
Thus we set $F_c  = \left(\begin{smallmatrix}1 & 0 \\ 0 & 1 \\ \sin\theta_{1}  & -\cos\theta_{1}  \\ \sin\theta_{2}  & -\cos\theta_{2}   \end{smallmatrix}\right)$ and we find $F_f =  \left(\begin{smallmatrix} -\sin\theta_{1} & -\sin\theta_{2}\\  \cos\theta_{1}  & \cos\theta_{2} \\  1 & 0\\    0 & 1 \end{smallmatrix}\right)$. We then set $F=(F_c\,|\, F_f)$. 

We do not include obstacles and thus  $b(\vx)\equiv 1$ and  
\begin{align*}
G &= H =F\operatorname{diag}([k\, k\, 1\, 1])F^\top=
&\scalebox{1.2}{$\left(\begin{smallmatrix} {\sin^2\theta_1}+{\sin^2\theta_2}+k & -\frac{\sin 2\theta_1}{2}-\frac{\sin 2\theta_2}{2} & (k-1)\sin\theta_1 & (k-1) \sin\theta_2 \\
 -\frac{\sin 2\theta_1}{2}-\frac{\sin 2\theta_2 }{2} & {\cos^2 \theta_1 }+{\cos^2 \theta_2 }+k & -(k-1)\cos \theta_1   & -(k-1)\cos \theta_2  \\ (k-1) \sin \theta_1   & -(k-1)\cos \theta_1   & k+1 & k\cos(\theta_1-\theta_2) \\ (k-1)\sin \theta_2   & -\cos \theta_2 \, k-1  & k\,\cos \theta_1-\theta_2  & k+1 \end{smallmatrix}\right)$}
 \end{align*}

\paragraph{The rolling coin or unicycle}

\begin{figure}\centering
  \subfloat[]{\label{fig:uniC}
           \centering
             \includegraphics[scale=1.2]{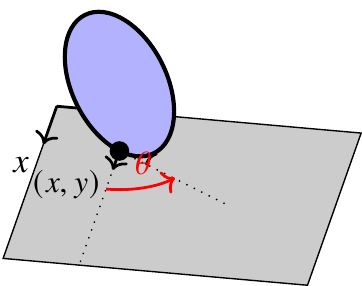}

      } 
\hspace{.3cm}\subfloat[]{
\includegraphics[scale=1.2]{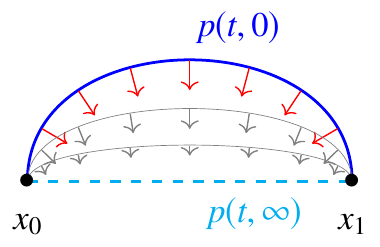}
\label{fig:meancurveflow}
       } 
 \caption{\small (a) A rolling coin or unicycle. is the side view. (b) In the mean-curvature flow, the curve $p(t,0)$ is continuously deformed in the direction of its normal, depicted by the red arrows. The final curve is a straight line. In general, the final curve is a length minimizing curve. is the corresponding angles} 
 \label{twofigs}
\end{figure}

The kinematics of a unicycle can be modeled as
\begin{equation}\label{unicycle}
\begin{pmatrix}\dot x\\\dot y\\\dot \theta\end{pmatrix}=\begin{pmatrix}\cos\theta\\\sin\theta\\0\end{pmatrix}u_1+\begin{pmatrix}0\\0\\1\end{pmatrix}u_2
\end{equation}
where $(x,y)$ is the position of the unicycle in the plane and $\theta$ is its orientation. Notice that there is only one non-holonomic constraints in this model and the constraint is the direction $\begin{pmatrix}-\sin\theta& \cos\theta &0\end{pmatrix}^\top$ which prevents moving sideways and hence prevents slipping. Equivalently, because the model~\eqref{unicycle} is affine in control, the free directions $F_f$ are simply the ones in~\eqref{unicycle}. Hence
\[
F(\vx)=\begin{pmatrix}
-\sin\theta&\cos\theta&0\\
\cos\theta&\sin\theta&0\\
0&0&1
\end{pmatrix},
\]
from which we obtain
\begin{equation*}
G(\vx)=H(\vx)=F \diag([k\ 1\ 1]) F^\top\\
=\begin{pmatrix}
\cos^2\theta+k\sin^2\theta&(1-k)\cos\theta\sin\theta&0\\
(1-k)\cos\theta\sin\theta&k\cos^2\theta+\sin^2\theta&0\\
0&0&1
\end{pmatrix}.
\end{equation*}

\subsection{Step 2: Initial sketch and solving the Geometric Heat Flow equation}

Our method proceeds with solving the following GHF equation:
\begin{equation}\label{HFE}
\frac{\partial}{\partial s}v_i(s,t)=\frac{\partial^2}{\partial t^2}v_i(s,t)+\sum_{j,k}\Gamma_{jk}^i\frac{\partial v_j}{\partial t}\frac{\partial v_k}{\partial t}\quad i=1,2,\dots,n
\end{equation}
where $\Gamma_{jk}^i$ are the Christoffel symbols introduced in \eqref{eq:defchristoffel} for the inner product defined in the previous subsection.  We impose  the boundary conditions
\[
v(s,0)=\vx_i,v(s,1)=\vx_f
\]
and a user defined initial condition, 
\[
v(0,t)=\vx(t)
\]
in order to find the solution. The initial curve $\vx(t)$ is an arbitrary curves satisfying the following 2 conditions:
\begin{enumerate}
\item It  satisfies the boundary conditions: $\vx(0)=\vx_i$ and $\vx(1)=\vx_f$;
\item It does not pass though any obstacles: $r(\vx(t))>0$ for all $t\in[0,1]$.
\end{enumerate}
An important point here is that $\vx(t)$ does not need to satisfy any holonomic or non-holonomic constraints; it can be simply a curve drawn from $\vx_i$ to $\vx_f$ without touching $\Omega$.

Notice that for each $s\geq 0$ fixed, the solution $v(s,\cdot)$ represent a curve connecting $\vx_i$ to $\vx_f$. As we explain below, as $s$ increases, $v(s,\cdot)$ is a curve that uses ``less and less of the constrained directions'', said precisely, $F_c^\top \frac{\partial}{\partial t} v(s,t)$ tends to zero.
We set $s_{\max}$ to be the simulation time for the PDE (in our examples, between 1 and 20) and  $$\vx_{sol}(\cdot)=v(s_{\max},\cdot).$$

\paragraph{Mean-curvature flows}\label{par:meanc} We now elaborate on the origin of Eq.~\eqref{HFE}: it is a type of curve-shortening flow~\cite{curvebook2001},  called a {\it mean-curvature flow} for a $1$-dimensional manifold (i.e. a curve) or {\it geometric heat flow}. For an  introduction to mean-curvature flows in arbitrary dimensions, see~\cite{colding2015mean}. For clarity of exposition, we present first the flow in two dimensional plane with the Euclidean inner product. We briefly mention steps that need to be taken for the general flow below. 

Consider a curve $p(t):[0,1] \to \R^2=(p_1(t),p_2(t))$, as depicted in Fig.~\ref{fig:meancurveflow}.  The scalar curvature~\cite{do1992riemannian} of $p$ at $p(t)$ is defined as $\kappa(p(t))=\|\ddot p\|$. Denote by $N_{p(t)}$ the unit normal vector pointing ``inward''.
The curvature of $p$ at $p(t)$ is then $\kappa(p(t)) N(p(t))$. 

The mean-curvature flow for this curve is defined as follows: consider a {\it family} of curves $p(t,s)$, $s \geq 0$, where for each $s_0$ fixed, $p(t,s_0):[0,1]\to\R^2$ is a curve joining $x_0$ to $x_1$, and $p(t,0)$ is the original curve.  Then the mean-curvature flow is the partial differential equation  $$\frac{\partial p}{\partial s} = \kappa(p(t,s)) N(p(t,s)).$$  Note that it is in fact a system of  two PDEs.
Looking at Fig.~\ref{fig:meancurveflow}, it is easy to conclude intuitively that $\lim_{s \to \infty} p(t,s)$ converges to a straight line between $x_0$ and $x_1$. This is also the {\it shortest path} between $x_0$ and $x_1$ for the usual Euclidean metric. This is no accident, and we can show that in general the solution of this PDE converges to a curve of minimal length.
For our purpose, we need to extend this idea in {\it two} directions: to $(i)$ curves in higher dimensions and $(ii)$ to a general Riemannian metric (or more precisely, inner product). One can show, after some extensive algebraic manipulations which we omit here, that the equivalent of the flow for a general curve in a  Riemannian manifold is  exactly the geometric heat flow presented in Eq.~\eqref{HFE}.

\subsection{Step 3: Extracting the controls}

The control can be directly computed:
\begin{equation}\label{control_extract}
\vu(t)=F_f^\dagger(\vx_{sol}(t))\dot \vx_{sol}(t)
\end{equation}
where $F_f^\dagger=(F_f^\top F_f)^{-1}F_f^\top$ is the pseudo-inverse of $F_f$. Notice that in the case $\vx_{sol}$ is admissible, that is, if $\dot \vx_{sol}(t)=F_fv(t)$ for some control $\vv$,
\[
\vu=F_f^\dagger \dot x_{sol}=(F_f^\top F_f)^{-1}F_f^\top F_f \vv=\vv
\]
Thus we have recovered the control and ideally the system should exactly follow the path $\vx_{sol}$. Notice that $F_fF_F^\dagger$ is a minimal square error projection onto the column space of $F_f$, the control extracted from \eqref{control_extract} will drive the system along a path that is close to $\vx_{sol}$, even if $\dot \vx_{sol}$ has small components in the constrained direction.

\subsection{On the implementation}
As mentioned earlier, the key of our method is to find an inner product matrix $G$ and then solve the GHF equation \eqref{HFE}. In our case, this is processed in MATLAB. To be explicit, once we have obtained $F_c$ from the constraints, we implement them as symbolic vectors in MATLAB and thus find $F_f(\vx)$. Subsequently, both $G(\vx)$ and $\frac{\partial G}{\partial \vx}$ can be derived symbolically and the symbolics are then replaced by state values and then stored in an $n\times n$ array \texttt{G} and an $n\times n\times n$ array \texttt{pG}, respectively. \texttt{pdepe} is then called with the boundary conditions and customized initial condition. In each iteration of solving the PDEs, the Christoffel symbols are computed from \texttt{G} and \texttt{pG} according to \eqref{eq:defchristoffel} and then stored in an $n\times n\times n$ array \texttt{Chris}. Notice that the \texttt{pdepe} solves PDEs of the general form
\[
c(s,t,x,\frac{\partial x}{\partial t})\frac{\partial x}{\partial s}=x^{-m}\frac{\partial}{\partial t}\left(t^m f(s,t,x,\frac{\partial x}{\partial t})\right)+s(s,t,x,\frac{\partial x}{\partial t})
\]
Compare it to~\eqref{HFE} we see that in our case we need to set\\
\texttt{c=ones(4,1);m=0,f=DxDt} and \texttt{s(i)=DxDt'*Chris(i,:,:)*DxDt}. Eventually the numerical solution of \texttt{pdepe} will be in the form of \texttt{sol(t,s,i)}, 

\subsection{Theoretical guarantee}

Set $\Delta(x) = \Span{\frac{\partial q_i}{\partial x}} \cap \Span{f_{c,j}}$.

We call the  constraints {\bf  satisfiable} if the distribution $\Delta$ satisfies the Lie algebraic rank condition (LARC). It is easy to see that it is a necessary condition for the existence of a trajectory joining arbitrary $\vx_i$ and $\vx_f$ while respecting the holonomic and non-holonomic constraints on the system. Under mild assumptions our method provides controls $\bar \vu(t)$ so that the solution $\vx^*(t)$ of $\dot \vx = \sum_i \bar \vu_i f_i$ by construction satisfies both the holonomic and non-holonomic constraints. In addition,

\begin{theorem}\label{main_thm}
Suppose $F(\vx)$ defined in \eqref{F(x)} is globally Lipschitz with constant $L$ and $\Vert F_c(\vx)\Vert=1$ for all $x\in \R^d$. Let $\bar E$ be the infimum of the energy functional
\[
E(\vu)=\int_0^1|\vu(t)|^2dt
\]
over the space of controls that the corresponding state trajectory satisfies both the holonomic and non-holonomic constraints. For any arbitrary $k\in \mathbb N,s>0$, define  $\vx$ to be the part $v(\cdot, s)$ of the solution of \eqref{HFE}, $\vu$ to be the control derived via \eqref{control_extract} and $\tilde \vx$ to be the solution of \eqref{eq:maindyn} generated by $\vu$ from $\tilde\vx(0)=\vx_i$. Then for any $\epsilon>0$, there exists $T=T(\epsilon,k)$ such that for all $s\geq T$,
\begin{enumerate}
\item $E(\vu)\leq \bar E+\epsilon$;
\item $|\tilde \vx(t)-\vx(t)|\leq \left(\sqrt{\frac{2t}{k}(\bar E+\epsilon)}\right)e^{L^2(\bar E+\epsilon)}$ for all $t\in[0,1]$. In particular, $|\tilde \vx(1)-\vx_f|\leq \left(\sqrt{\frac{2}{k}(\bar E+\epsilon)}\right)e^{L^2(\bar E+\epsilon)}$.
\end{enumerate}
\end{theorem}

\section{Case study} 

\paragraph{Articulated arm}

\begin{figure}\centering
  \subfloat[]{
           \centering
     \includegraphics[scale=1.2]{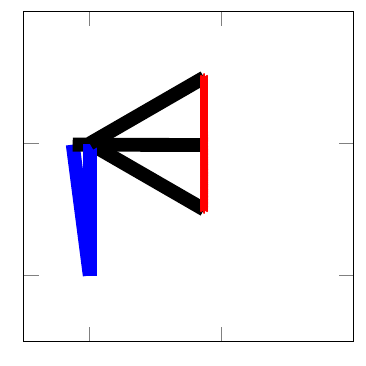}
} 
\hspace{.3cm}\subfloat[]{
       \includegraphics[scale=1.2]{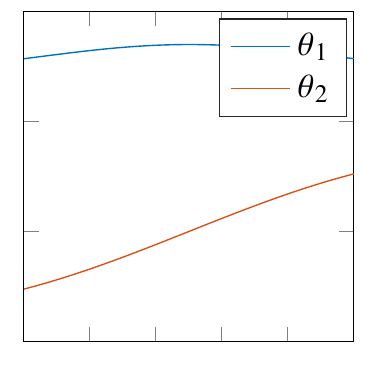}
}

  \subfloat[]{
           \centering
      \includegraphics[scale=1.2]{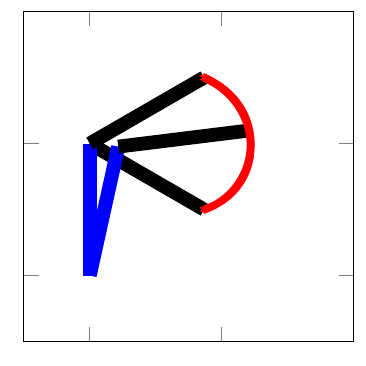}
} 
\hspace{.3cm}\subfloat[]{
       \includegraphics[scale=1.2]{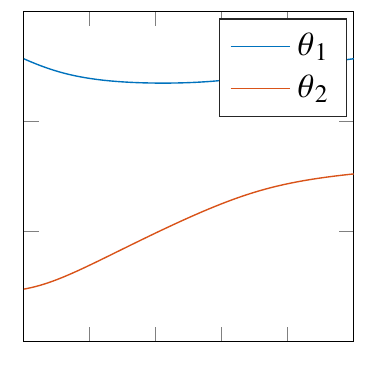}
} 
 \caption{\small Vertical motion (a) and circular  motion (c) of the two-links articulated arm. The links are in blue and black. The trajectory of the tip is marked in red. We draw the initial and final state and an intermediate state. The joint angles are given in (b) and (d) respectively. 
  \label{two_links} } 

\end{figure}

We first study the 2R robot introduced earlier. Our goal is to plan the motion of the tip of the arm, from an initial state $\vx(0)=\vx_i=(\sqrt{2}/2,1-\sqrt{2}/2,\pi/2,-\pi/4)$, where we recall that the coordinates are $(x,y,\theta_1,\theta_2)$, to a final state  $\vx(1)=\vx_f=(\sqrt{2}/2,1+\sqrt{2}/2,\pi/2,\pi/4)$. We furthermore require the motion to follow a {\it straight line} given by $x=constant$. The resulting motion planning problem thus contains,
in addition to the two holonomic constraints relating the tip of the arm to the angles given in Eq.~\eqref{eq:twoarmconst}, the  constraint $q_3(\vx)=x-x_i=0$ and the corresponding constrained direction is $\frac{\partial q_3}{\partial \vx}=(1,0,0,0)^\top$. Given these constraints, we implement the three steps of the method outlined above show the results in Fig.~\ref{two_links}. We then replaced the constraint of vertical motion by asking that the tip follows an arc of a circle. The corresponding holonomic constraint is $q_4(\vx)=(x-x_c)^2+(y-y_c)^2 - r=0$ for some constants $x_c,y_c,r$. The differential of this constraint is easily evaluated. We show in Fig.~\ref{two_link_circular}  the result obtained. Note that this illustrate the use of our method to solve {\it inverse kinematic problems} numerically.

\begin{figure}[ht!]\centering
    \subfloat[Initial sketch in $(x,y,\theta)$-space. ] {
 \includegraphics[scale=.95]{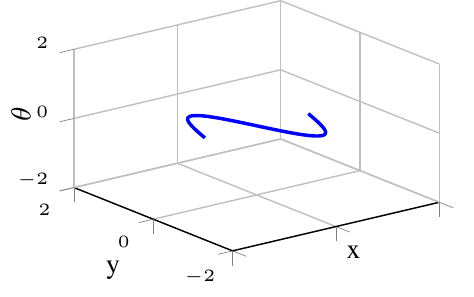}\label{3D_initial}}  
\hspace{.3cm}
\subfloat[Projection of the initial sketch in $(x,y)$- plane, with snapshots of the corresponding position of the unicycle. Note that $\theta=0$ for each snapshot.]{
 \centering \includegraphics[]{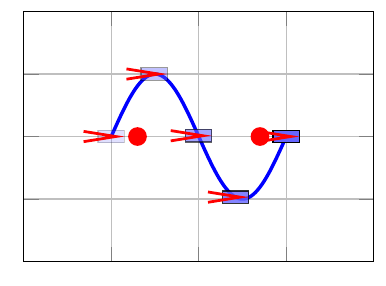} \label{2D_initial}}
  
  \subfloat[Solution $\vx(s,t)$ for $s=0.0001$]{
 \includegraphics[]{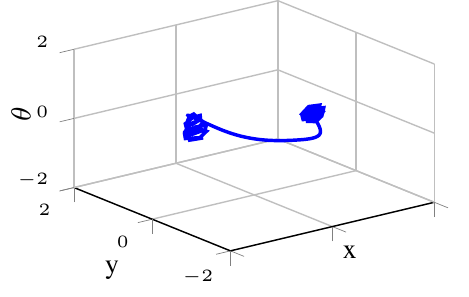}\label{3D_0.0001}} 
 \hspace{.58cm}
 \subfloat[Plane view at $s=0.0001$]{
  \centering \includegraphics[]{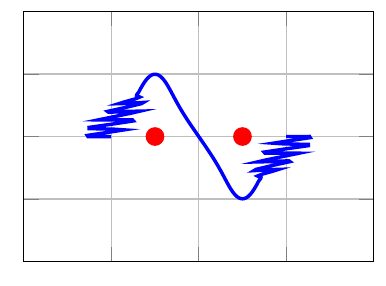}  \label{2D_0.0001}}       
 
 \subfloat[Solution $\vx(s,t)$ for $s=0.01$.]{
 \includegraphics[]{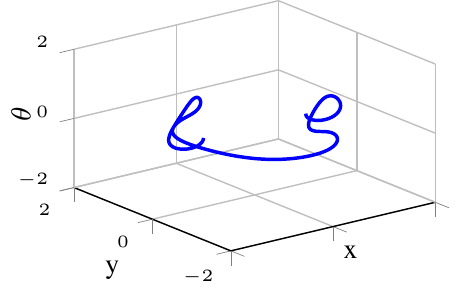}\label{3D_0.01}} 
  \hspace{.58cm}
  \subfloat[Solution $\vx(s,t)$ for $s=0.01$]{
  \centering \includegraphics[]{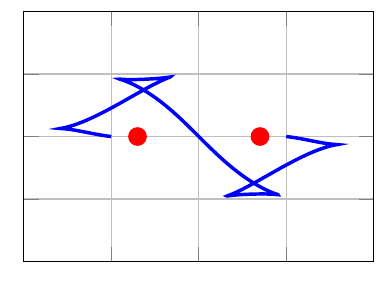}  \label{2D_0.01}}    
  
  \subfloat[3D view of the final curve]{
 \includegraphics[]{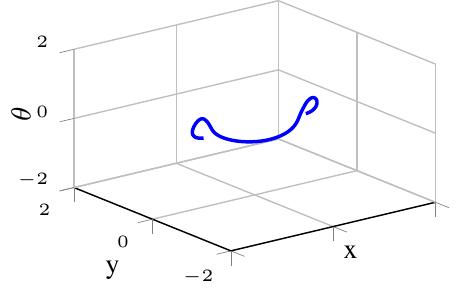}\label{3D_final}} 
  \hspace{.58cm}
  \subfloat[Plan view of the final curve]{
  \centering \includegraphics[]{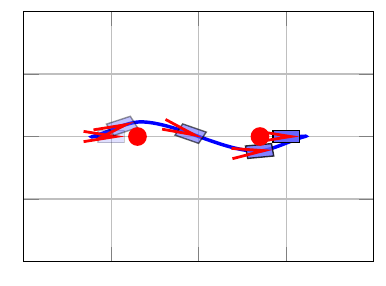} 
  \label{2D_final}}                          
   
   \caption{\small Path planning for a unicycle avoiding two point obstacles. The red dots are the two obstacles, the blue curves are the solution of GHF equations at different $s$. In the plan views of initial curve and final curve, unicycle positions are marked along the curve, with its orientation labeled with red arrows }\label{unicycle_solution}
\end{figure}

\paragraph{Unicycle}
Consider the unicyle described above with coordinates $(x,y,\theta)$. We desire to transfer the unicycle from $(x(0),y(0),\theta(0))=(-1,0,0)$ to  $(x(1),y(1),\theta(1))=(1,0,0)$ without slip (a non-holonomic constraint). In addition, there are two point obstacles located at $(-0.7,0),(0.7,0)$ which the unicycle should avoid in the xy-plane. Provided these constraints, we first build an inner product $G(x)$ as described earlier. We then provide  an arbitrary curve connecting $(-1,0,0)$ and $(1,0,0)$ and avoiding the obstacles--we called this curve the initial sketch. We opted simply for a sinusoidal curve in xy-plane and kept  $\theta\equiv 0$, as shown in Fig.~\ref{3D_initial}. As observed in Fig.~\ref{2D_initial}, the unicycle certainly cannot follow this curve, as the motion direction is not aligned with the unicycle orientation or, in other words, the non-slip constraint is not met. 

Recall that the solution of GHF equations~\eqref{HFE} is a curve connecting the initial and final states when $s$ fixed. In Figs.~\ref{3D_0.0001} to~\ref{3D_final}, we show the gradual  deformation of the curve in configuration space as $s$ increases. In the final step  $s=4$, the curve becomes almost admissible and we see that the unicycle can basically follow such trajectory to reach its final state. It is worth noticing that because the obstacles are very close to the initial and final states, the unicycle has to move backward first in order to have more room to maneuver around said obstacles. Similarly, it overshoots the second obstacles before backing up and parking at its final destination. 

\paragraph{Car}

\begin{figure}
        \centering
        \subfloat[Car modelled by Eq.~\eqref{boxcar_dynamics}. The red arrow is used to indicate the front of the car.
        ]{  
   \includegraphics[scale=1.3]{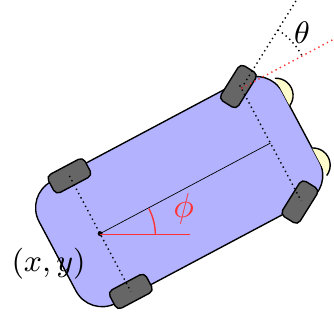}  \label{fig:boxcar} 
        
        }\hspace{.58cm}
   \subfloat[Initial sketch. The car rotates 180 degrees with its center of mass following the black curve with slipping.]{

        \includegraphics[scale=.83]{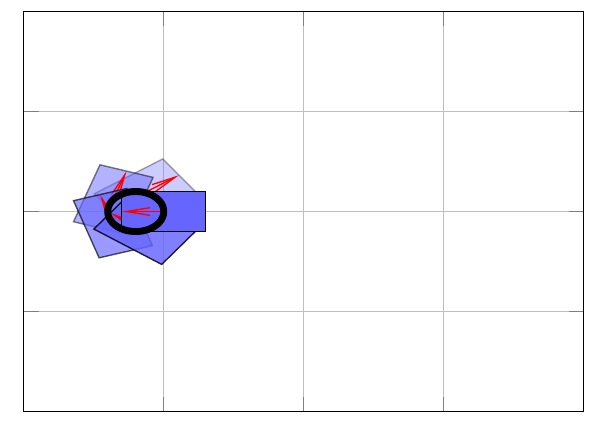} \label{car_180_initial}}
        
  \subfloat[{3 points turning when no spatial constraints}]{
           \centering
      \includegraphics[scale=1.3]{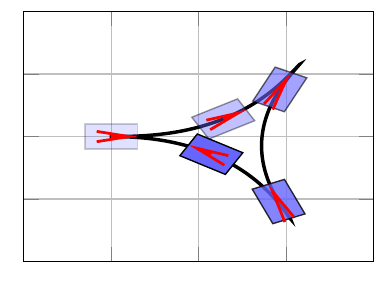} \label{car_180_free}}  
\hspace{.41cm}\subfloat[5 points turning between walls]{
       \includegraphics[scale=1.3]{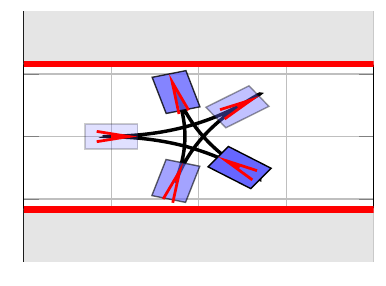} \label{car_180_street}}
 \caption{Car $180^\circ$ turn experiment.} 
 \label{car_180}
\end{figure}

We now illustrate our method for planning the motion of a car with position $(x,y) \in \R^2$, body orientation $\phi$ and wheel angle $\theta$. A top view of car is illustrated in Fig.~\ref{fig:boxcar} and the equations of motion equation are:
\begin{equation}\label{boxcar_dynamics}
\begin{pmatrix}\dot x\\\dot y\\\dot \theta\\ \dot \phi\end{pmatrix}
=u_1\begin{pmatrix}\cos\phi\\\sin\phi\\0\\\frac{1}{d}\sin\theta\end{pmatrix}+u_2\begin{pmatrix}0\\0\\1\\0\end{pmatrix},
\end{equation}
where $u_1$ is the throttle input, $u_2$ is the steering input and $d$ is the distance between front wheels axis and rear wheels axis. We have studied this example in our paper~\cite{Belabbas2017NewMF}, and we refer the reader to this paper for an explicit derivation of the corresponding $G(\vx)$.

Our first experiment is  a $180^\circ$  turn. Our initial sketch for this motion is illustrated in Fig.~\ref{car_180_initial}. It is clear that, unless equipped with omniwheels or $d=0$,  the car cannot perform the motion illustrated. Interestingly, Motionsketch deforms this curve into the well-known  3-points turn path illustrated in Fig.~\ref{car_180_free}. This corresponds to the most efficient way of $180^\circ$ turning of a car in practice, assuming there are no any other spatial obstacles. 

If in addition, we impose add parallel curbs, which are encoded in the barrier function $b(\vx)$ as described earlier,  the constrained space the car can move in results in additional back-and-forth.  The narrower the street, the more back-and-forth are needed. We provide additional examples in the webpage\footnote{\url{https://publish.illinois.edu/belabbas/motion-planning/}}.

\begin{figure}
        \centering
   \subfloat[Initial sketch. Note that the constraints are not met.]{
        \centering
        \includegraphics[scale=.7]{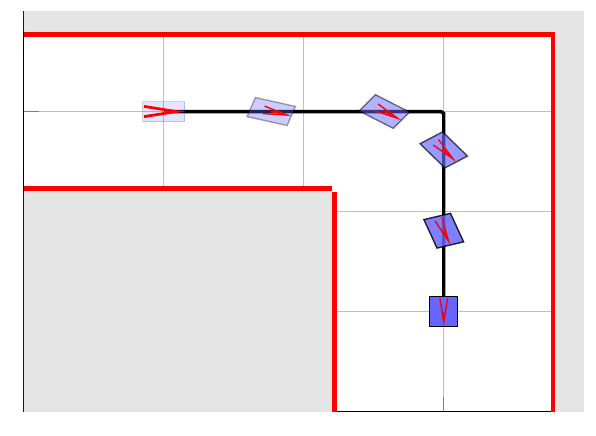} \label{car_90_initial}}
  \subfloat[Turn in a wide street corner]{
           \centering
      \includegraphics[scale=1]{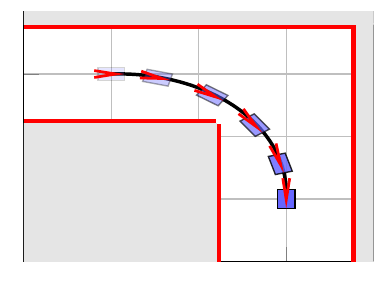} \label{car_90_wide}}  
\hspace{-.0cm}\subfloat[Back-forth behavior at narrow street corner]{
       \includegraphics[scale=1]{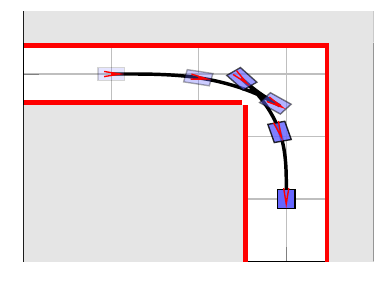}  \label{car_90_narrow}}
\caption{\small Car street corner turn experiment}
\label{car_90}
\end{figure}

We conclude with  the case of a car turning in a narrow street. The initial curve is simply an L-shaped curve in xy-plane with  $\phi$ linear with respect to $t$ and $\theta\equiv 0$, as illustrated in Fig.~\ref{car_90_initial}. With the curbs modeled as obstacles, our method generates the relatively ``optimal'' path for this corner turn. Interestingly enough, the car is able to perform the turn in one shot if the street is relatively wide as shown in Fig.~\ref{car_90_wide}, or may need extra maneuvering  if the street is narrow, as shown in Fig.~\ref{car_90_wide}. We emphasize that  both simulation are performed with the same initial curve provided in Fig.~\ref{car_90_initial}. The only difference is the street width. Whether one shot or two is automatically determined by our method without any further specification. 

Finally, we note that in addition to the curb of the streets which are modeled as obstacles in the xy-plane, we also put limits on the steering angle $\theta$ as an obstacle for the $\theta$ variable.

\paragraph{Multi-vehicle path planning}
We show that multiple vehicles can be path planned simultaneously using our methods. In the first simulation two unicycles are initially at states $(0,1,0),(0,-1,0)$; that is, parked at xy-coordinate $(0,1),(0,-1)$ while both facing east. The task is to swap the position of the two unicycles. The initial sketch is a circle passing through the two unicycles -- clearly these two paths are infeasible since the orientation vectors of the unicycles are not tangent to the paths. After running our algorithm, the two initial sketch of paths deform into the two V-shaped paths and now the two unicycles are able to perform the swap of positions along such paths {\it while avoiding collisions}.
While readers might think the previous example has no major difference compared with path planning for single vehicle and hence less challenging, the next example is more interesting and shows the power of our algorithm in multi-vehicle path planning. In this case one unicycle is supposed to move from $(-1,0,\pi/2)$ to $(1,0,\pi/2)$ while the other one is supposed to move from $(0,-1,0)$ to $(0,1,0)$.

\begin{figure}\centering
\noindent\subfloat[]{
           \centering
     \includegraphics[scale=1.2]{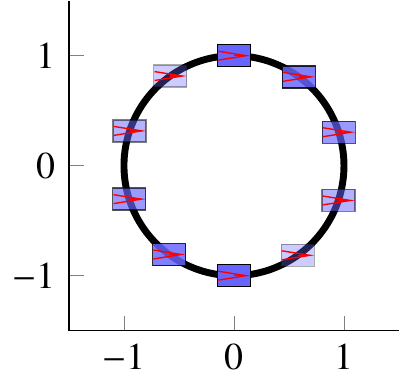} }
\hspace{.3cm}\subfloat[]{
       \includegraphics[scale=1.2]{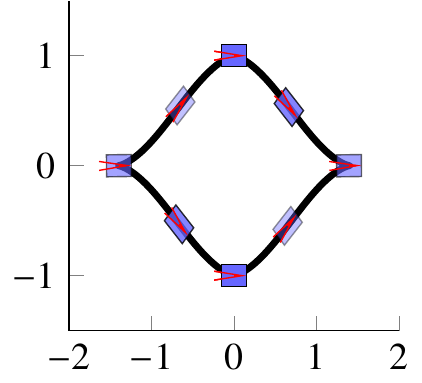}}  
 \subfloat[]{
           \centering
     \includegraphics[scale=1.2]{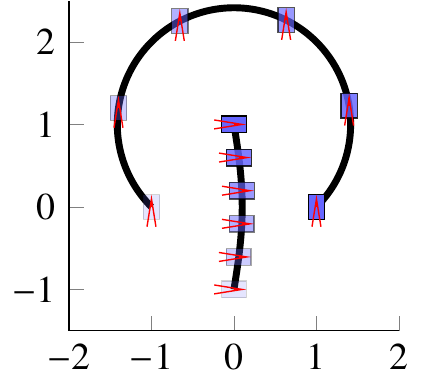} } 
\hspace{.3cm}\subfloat[]{
       \includegraphics[scale=1.2]{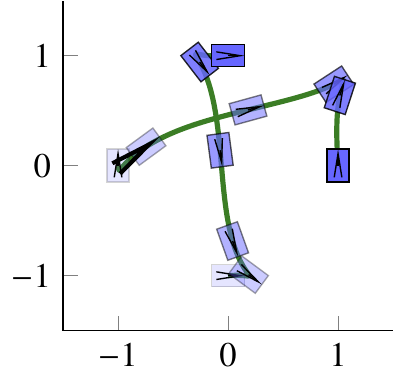} } 
 \caption{\small Multi-vehicles motion planning with collision avoidance.
  \label{two_cars} } 

\end{figure}

\section{Summary and discussion}

We have provided in this paper a guide to the implementation of the method we termed MotionSketch for solving motion planning problems. We have illustrated the use of the method on examples with holonomic, non-holonomic and obstacle constraints, and have demonstrated that the method yields good practical results.

The  salient points  of the method were that it encodes all the constraints into a Riemannian inner product, and that it requires an initial sketch of the curve joining a desired final state to an initial state. This curve however does not need to meet the holonomic and non-holonomic constraints and is thus often easily obtained. In fact, if the space is convex, a straight line joining the two states most often meets the constraints.

Amongst the problems that are also readily solved using MotionSketch, but that we did not show here, we mention multi-vehicle motion planning with collision avoidance. For example, think of having to plan the trajectory of two non-holonomic cars with the constraints that they should avoid each other. This can be done using our method as follows: denote by $(x_i,y_i,\theta_i,\phi_i)$ the coordinates describing the state of car $i$, and by $G_i\in \R^{4 \times 4}$ the corresponding Riemannian inner products modeling the constraints for each car (e.g. max turning angle as am obstacle in $\theta$, curbs, etc.). In order to model the two vehicles scenario, we first consider the cartesian product of the coordinates with metric $\bar G \in \R^{8 \times 8}$ a block diagonal matrix with blocks $G_i$. In order to avoid collisions between the cars, it suffices to place an obstacle around the ``diagonal'' subspace $x_1=x_2$ and $y_1=y_2$. As we have seen earlier, adding obstacles to a metric only requires multiplying by a barrier function, hence we can set $G(x) = b(x) \bar G(x)$. This procedure  generalizes in a straightforward way to the case of more than two vehicles.

\paragraph{On the computational complexity of solving the GHF}

The numerically intensive part of the method lies in solving the geometric heat flow, which is a system of parabolic partial differential equations. We point out that solving such a PDE can be done rather efficiently, owing to the fact that the complexity scales polynomially with the dimension, and not exponentially, and the fact that there exists parallel algorithms to do so. 

To elaborate on the first point, the main reason why  the PDE we use scales well is that the {\it domain} of its solution has a {\it constant} dimension of two. For most PDEs encountered in engineering, such as the heat equation, or the Hamilton-Jacobi-Bellman equation, the dimension of the problem affects the  dimension of the {\it domain} of the solution seeked, whereas is our case, it affects the dimension of the image of the solution. A linear increase in the dimension of the domain yields what is often referred to as the {\it curse of dimensionality}, as the number of interpolation points needed to represent a function on a domain of dimension $n$ grows exponentially with $n$. Note however that the domain of our PDE is {\it always} two-dimensional, but the dimension of the image increases linearly, the number of interpolation points grows {\it linearly} with the dimension.  Hence our PDE does not suffer from the curse of dimensionality and thus scales well to higher-dimensional problems. We refer to, e.g.,~\cite{} for a more detailed discussion on the complexity of solving such PDEs. In practice, using MATLAB on a common laptop computer with non-optimized code (in particular, MATLAB does not solve such PDEs using multiple cores), the computation time was of the order of seconds to minutes, depending on the complexity of the problem. Per our discussion above, we believe however that there is ample room for improvement on this front.

\end{document}